\documentclass{llncs}
\usepackage{amssymb,amscd, epsfig}
\title{
Non-standard discretization 
of biological models} 
\author{  
Andrew Hone and Kim Towler   
} 
\authorrunning{A. Hone and K. Towler}
\institute{School of Mathematics,
Statistics \&
Actuarial Science, \\ 
University of Kent, Canterbury CT2 7NF, UK.
\\
\email{anwh@kent.ac.uk}
} 

\begin{document}

\newcommand{\haf}{{\hat{f}}}
\newcommand{\beq}{\begin{equation}}
\newcommand{\eeq}{\end{equation}}
\newcommand{\bea}{\begin{eqnarray}}
\newcommand{\eea}{\end{eqnarray}}
\newcommand\la{{\lambda}}
\newcommand\ka{{\kappa}}
\newcommand\al{{\alpha}}
\newcommand\be{{\beta}}
\newcommand\de{{\delta}}
\newcommand\si{{\sigma}}
\newcommand\lax{{\bf L}}
\newcommand\mma{{\bf M}}
\newcommand\ctop{{\mathcal{T}}}
\newcommand\hop{{\mathcal{H}}}
\newcommand\ep{{\epsilon}}
\newcommand\T{{\tau}}
\newcommand\om{{\omega}}
\newcommand\ga{{\gamma}}
\newcommand\zbar{{\overline{z}}}
\newcommand\tal{{\hat{\alpha}}}
\newcommand\tbe{{\hat{\beta}}}
\newcommand\mez{{\frac{1}{2}}}

\newcommand\rd{\mathrm{d}}

\newcommand{\bear}{\begin{array}}
\newcommand{\eear}{\end{array}}

\def\endpf{\begin{flushright}$\square$\end{flushright}}
\newtheorem{prop}[theorem]{Proposition}
\newtheorem{lem}[theorem]{Lemma}
\newtheorem{cor}[theorem]{Corollary}
\newtheorem{defi}[theorem]{Definition}

\newenvironment{prf}{\trivlist \item [\hskip
\labelsep {\bf Proof:}]\ignorespaces}{\qed \endtrivlist}

\newcommand\C{\mathcal{C}}
\newcommand\F{\mathcal{F}}
\newcommand\G{\mathcal{G}}

\newcommand\PP{\mathcal{P}}

\maketitle
\begin{abstract}
We consider certain types of discretization schemes 
for differential equations with quadratic nonlinearities, which 
were introduced by Kahan, and considered in a broader  
setting by Mickens. 
These methods have the property that they preserve 
important structural features of the original systems, 
such as the behaviour of solutions near to fixed points, 
and also, where appropriate (e.g. for certain mechanical systems),  
the property of being volume-preserving, or preserving a symplectic/Poisson structure.  
Here we focus on the application of Kahan's method to models of biological systems, in particular 
to reaction kinetics governed by the Law of Mass Action, 
and present a general approach to birational discretization,  which is applied 
to population dynamics of Lotka-Volterra type.  
\end{abstract}

\section{Introduction}

In 1993 Kahan gave a set of lectures entitled ``Unconventional numerical methods for trajectory calculations,'' 
in which he proposed a method for discretizing a set of differential equations 
\beq \label{ode}
\dot{\bf x}={\bf f}({\bf x}), 
\eeq 
where all the components of the vector field ${\bf f}$ are polynomial functions of 
degree at most two in the components $x_1,x_2,\ldots ,x_N$ of the vector ${\bf x}$, 
and the dot denotes the time derivative $\rd /\rd t$. 
Kahan's method consists of replacing the left hand side of (\ref{ode}) by the standard 
forward difference, while on the right hand side quadratic, linear and constant terms are 
replaced according to a symmetric rule, so that overall 
the method is specified as follows: 
\beq\label{kmet} 
\dot{x}_i\to \frac{\tilde{x}_i-x_i}{h}, \quad 
x_ix_j\to\frac{x_i\tilde{x}_j+ \tilde{x}_ix_j}{2}, \quad 
x_i\to\frac{x_i+ \tilde{x}_i}{2}, \quad c\to c.  
\eeq    
Above and throughout we use the tilde to denote 
the finite difference approximation to a dependent variable 
shifted by a time step $h$,  
i.e. $x_i(t+h)\approx \tilde{x}_i$. 

The replacements (\ref{kmet}) result in a difference 
equation of the form 
\beq \label{disc} 
\frac{\tilde{{\bf x}}-{\bf x}}{h} = {\bf Q} ({\bf x}, \tilde{{\bf x}}), 
\eeq 
where the right hand side is a vector  function of degree two. 
Thus it appears to be an implicit scheme, in the sense that (\ref{disc}) 
defines $\tilde{{\bf x}}$ 
implicitly as a function of 
${\bf x}$. However, the fact that 
the formulae (\ref{kmet}) are linear in each of the variables 
$\tilde{x}_i$ and $x_j$ 
means that (\ref{disc}) can be solved explicitly to find  
$\tilde{{\bf x}}$ as a rational function of ${\bf x}$, 
and vice versa, yielding a \textit{birational} map 
$\varphi:  {\bf x}\mapsto \tilde{{\bf x}}$. 
As shown in \cite{KHLI}, the map can be written explicitly as 
\beq \label{phi} 
\varphi: \quad \tilde{{\bf x}} = {\bf x} + h 
\left({\bf I}-\frac{h}{2}{\bf f}'({\bf x})\right)^{-1}{\bf f}({\bf x}), 
\eeq 
where ${\bf I}$ denotes the $N\times N$ identity matrix, and ${\bf f}'$ is the Jacobian of ${\bf f}$, 
while the inverse is 
\beq \label{phinv} 
\varphi^{-1}: \quad {{\bf x}} = \tilde{{\bf x}} - h \left({\bf I}+\frac{h}{2}{\bf f}'(\tilde{{\bf x}})\right)^{-1}{\bf f}(\tilde{{\bf x}}). 
\eeq 
The above method is second-order, but Kahan  and  Li showed
how it can be used with suitable composition schemes 
to generate methods of higher order \cite{KHLI,kahan}. 
  
Roegers proved that (in contrast to Euler's method) Kahan's 
method preserves the local stability of steady states of (\ref{ode}). 
Indeed, the steady states ${\bf x}^*$ of the differential system,  
satisfying ${\bf f}({\bf x}^*)=0$, coincide with those of 
(\ref{phi}), i.e. the solutions of $\varphi ({\bf x}^*) = {\bf x}^*$, 
and taking the derivative of $\varphi$ at such an ${\bf x}^*$ gives 
$$ 
\varphi '({\bf x}^*) = {\bf I} + h\left({\bf I}-\frac{h}{2}{\bf f}'({\bf x}^*)\right)^{-1}{\bf f}'({\bf x}^*) .  
$$ 
Hence to each eigenvalue $\la$ of ${\bf f}'$ at ${\bf x}^*$ there corresponds 
an eigenvalue $\mu (h)$ of $\varphi '$ at ${\bf x}^*$, with the same eigenvector, where 
\beq \label{mu} 
\mu (h)= \frac{1+\frac{h\la}{2}}{1-\frac{h\la}{2}}. 
\eeq  
The above transformation sends the region  Re$\, \la <0$ 
to  $|\mu (h)|<1$, which identifies asymptotically stable directions at steady states of 
(\ref{ode}) with those for (\ref{phi}), and similarly for unstable directions (Re$\, \la >0$ is sent 
to $|\mu (h)|>1$). These local stability properties go some way towards explaining why 
Kahan's method seems to preserve global structural features of solutions of (\ref{ode}). 

In the context of Hamiltonian mechanics, Hirota and Kimura rediscovered 
Kahan's prescription (\ref{kmet}) as a new method to discretize Euler's equations 
for a top spinning about a fixed point \cite{HK1}. This stimulated further  interest in 
the method, and as a result many new 
completely integrable symplectic/Poisson maps were found (see \cite{hp}, for instance). 
Despite an extensive survey of algebraically completely integrable discrete systems obtained   
via Kahan's method \cite{PPS2}, the general conditions 
under which a quadratic Hamiltonian vector field (\ref{ode}) produces a 
map (\ref{phi}) that preserves a symplectic (or Poisson) structure, as
well as one or more first integrals, are still unknown. 
Nevertheless, considerable progress was made recently by Celledoni et al. \cite{celledoni}, 
who showed that Kahan's method is perhaps not as ``unconventional'' as originally thought, 
since it coincides with the Runge-Kutta method 
$$ 
 \frac{\tilde{{\bf x}}-{\bf x}}{h} = -\frac{1}{2}{\bf f}({\bf x}) + 2 {\bf f} \left( \frac{{\bf x}+ \tilde{{\bf x}}}{2}\right) 
-\frac{1}{2}{\bf f}(\tilde{{\bf x}}) 
$$ 
applied to quadratic vector fields ${\bf f}$. Moreover, if (\ref{ode}) is a Hamiltonian system 
with a constant Poisson structure and a cubic Hamiltonian function, then the corresponding 
map (\ref{phi}) has  a rational first integral and preserves a volume form (Propositions 
3,4 and 5 in \cite{celledoni}).


Rather than applications in mechanics, in this paper we focus on some applications of Kahan's method to 
biological models. In the next section we discuss how this method is well suited to modelling reaction kinetics. 
The basic enzyme reaction is used as an example to illustrate the method, and we see how the discretization 
reproduces the transient behaviour inherent in this system. 

Mickens proposed a broad non-standard approach to preserving structural 
features of differential equations under discretization \cite{mickensbook}, which has been applied 
to many different problems (as reviewed by Patidar in \cite{patidar}).     
In the third section we consider population dynamics, and more specifically the Lotka-Volterra 
model for a predator-prey interaction. This was one of the examples originally treated by Kahan 
in applying his symmetric method, but Mickens found an asymmetric discrete Lotka-Volterra 
system with the same qualitative features \cite{mickens}. Here we give the details of a general method, 
first sketched in \cite{hone},  
to obtain non-standard discretizations by requiring that the resulting  maps should be \textit{birational}.  
As an application of the method, a classification of birational discrete Lotka-Volterra systems is derived, 
and this is then used to reproduce some results of Roeger on symplectic Lotka-Volterra maps \cite{roeger2}.  
 
The final section is devoted to some conclusions. 

\section{Discrete reaction kinetics: the basic enzyme reaction} 

Systems of differential equations of the form (\ref{ode}), with ${\bf f}$ being a quadratic vector field, 
are ubiquitous in chemistry, and in biochemistry in particular. They arise  
immediately from reaction kinetics  involving reactions of the form $A+B\mathop{\rightarrow}C$, 
$A\mathop{\rightarrow}B+C$, or $A+B\mathop{\rightarrow}C+D$, where $A,B,C,D$ represent different 
molecular species.
In that case, the Law of Mass Action implies that 
the rate of change of concentration of each reactant 
is given by a sum of linear and quadratric terms in the concentrations. 
Processes involving collisions of three or more molecules are statistically rare, and 
although trimolecular reactions may be observed empirically, it is usually understood that 
in practice they are mediated by dimolecular reactions (which may take place very rapidly). 
Thus quadratic nonlinearities are the norm in reaction kinetics.  

Kahan's discretization method is well-suited to  reaction kinetics models, 
where the variables $x_i$, $x_j$ in  (\ref{kmet}) would represent concentrations of different 
chemical species. From a practical point of view, the only potential difficuly is with inverting 
the matrix ${\bf I}-\frac{h}{2}{\bf f}'({\bf x})$ on the right hand side of  (\ref{phi}), 
which becomes unfeasible to do algebraically when $N$ (the number of species) is large. 
However, viewed numerically for a given value of ${\bf x}$ 
and a small $h$, this matrix is a small perturbation of the identity, and its inverse can be expanded as 
a geometric series, 
$$ 
\Big({\bf I}-\frac{h}{2}{\bf f}'\Big)^{-1} = 
{\bf I} + \frac{h}{2} {\bf f}' +  \frac{h^2}{4} ({\bf f}')^2 + \frac{h^3}{8} ({\bf f}')^3 + \ldots , 
$$ 
and then truncated at a suitable power of $h$ if necessary. 
Moreover, if each species in the reaction network is only coupled to a small number 
of others then the matrix will be sparse, and there are efficient algorithms for inverting 
matrices of  this  kind.

To illustrate how Kahan's method works in reaction kinetics, we 
consider the basic enzyme reaction, 
which is given by
$$ S+E \mathop{\rightleftarrows}^{k_1}_{k_{-1}} C,
\qquad C \mathop{\longrightarrow}^{k_2} E+P,$$
where $S$ is the substrate, $E$ is the enzyme, $C$
is a combined substrate-enzyme complex, $P$ is the product and
$k_{-1},k_1,k_2$ are rate constants. Our presentation follows that of chapter 6 in \cite{murray2} 
very closely. 

From the above reaction scheme, the Law of Mass Action yields a system 
of four differential equations which describe the rate of change of 
concentration of each reactant, as follows:   
\beq \label{der}
\begin{array}{rclrcl} 
\dot{s} & = & -k_1 es + k_{-1}c, \qquad \qquad  & \dot{e} & = & -k_1 es + (k_{-1}+k_2)c, \\ 
\dot{c} & = & k_1 es - (k_{-1}+k_2)c, & \dot{p} & = & k_2 c. 
\end{array}  
\eeq 
The lower case letters $s,e,c,p$ are used to denote the concentrations 
of the reactants $S,E,C,P$, respectively. For an enzyme reaction, we can assume 
that initially none of the enzyme has bound with the substrate, and no product has yet formed,  
so that the initial conditions take the form  
\beq \label{inits} 
 s(0)=s_0, \qquad e(0)=e_0, \qquad c(0)= 0 = p(0) . 
\eeq

The right hand side of the system (\ref{der}) has degree two overall, which means 
that Kahan's method (\ref{kmet}) can be applied immediately, to produce the   
discretization 
\beq\label{kder} 
\begin{array}{rcl} 
(\tilde{s}-s)/h & = & - k_1(e\tilde{s}+\tilde{e}s)/2 +k_{-1}(c+\tilde{c})/2,  \\   
(\tilde{e}-e)/h & = & - k_1(e\tilde{s}+\tilde{e}s)/2 + (k_{-1}+k_2)(c+\tilde{c})/2, \\ 
(\tilde{c}-c)/h & = & k_1(e\tilde{s}+\tilde{e}s)/2 - (k_{-1}+k_2)(c+\tilde{c})/2, \\  
(\tilde{p}-p)/h & = &  k_2(c+\tilde{c})/2. 
\end{array}
\eeq 
However, rather than go ahead and solve the discrete system explicitly for 
$\tilde{s}, \tilde{e},\tilde{c},\tilde{p}$, we will use a special feature 
of the enzyme reaction in order to simplify the problem.  

The essential feature of an enzyme is that it is a catalyst, so that it is not changed overall 
by the reactions it is involved in. This means that the total amount of enzyme 
in the system (free enzyme in the form of $E$, and enzyme bound into the complex $C$) 
must be conserved with time. From (\ref{der}) this can be seen directly in the form 
of the conservation law 
$$ 
\dot{e} + \dot{c} =0 \qquad  
\Rightarrow 
\qquad e(t) + c(t) = \mathrm{constant} = e_0, \qquad 
\forall t\geq 0,   
$$ 
 where  the initial conditions (\ref{inits}) are used to specify the value of 
the first integral $e+c$. 

\begin{figure}\label{pxyz} 
\begin{center}$ 
\begin{array}{cc} 
{\includegraphics[width=2.2in]{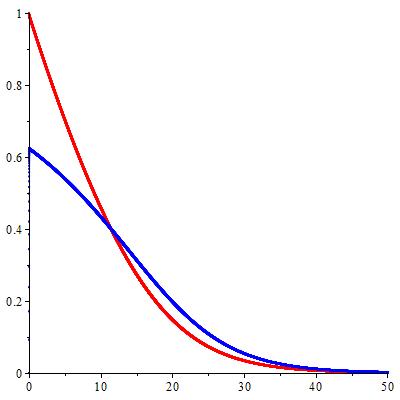}} & 
{\includegraphics[width=2.2in]{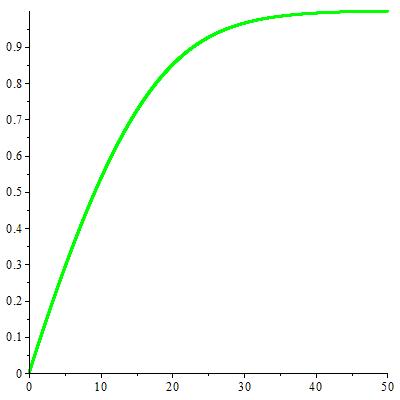}} 
\end{array}$ 
\end{center}
\caption[ ]{Plots of $x$ and $y$  against $t$ (left) and $z$ (right) against $t$ from numerical integration of (\ref{diml}) using Kahan's method 
with $\nu=0.6$, $\mu=0.5$, $\varepsilon =10^{-2}$, $h=10^{-3}$.}    
\end{figure}  

The fact that the linear function $e+c$ is conserved means that $e$ (or $c$) can be eliminated 
from the system (\ref{der}), so that the problem is reduced to solving a three-dimensional system. 
As noted in \cite{celledoni}, the properties of Kahan's method mean that it preserves linear 
first integrals, and so, as is clear from the above equations, the discrete system (\ref{kder}) 
has the same first integral: 
$$ 
\tilde{e} + \tilde{c} = e+c. 
$$ 
Thus we can also set $e=e_0-c$ to eliminate $e$ from (\ref{kder}), 
and then solve the resulting three-dimensional discrete system. 
Rather than doing so immediately, it is helpful to rescale 
all the variables in (\ref{der}), so that (after eliminating $e$) everything 
is written in terms of dimensionless quantities $x,y,z$ and $\tau$, where 
$$ 
x=\frac{s}{s_0}, \qquad y = \frac{c}{e_0}, \qquad z = \frac{p}{s_0}, \qquad 
\tau = k_1 e_0 t.    
 $$ 
Then the dimensionless system is 
\beq \label{diml}
\begin{array}{rcl}
\dot{x} & = & -x + \mu y + xy, \\ 
\varepsilon \dot{y}  & = &  x- \nu y - xy, \\ 
\dot{z} & = & (\nu - \mu )y, 
\end{array}
\eeq 
where now the dot denotes $\rd /\rd \tau$, and the dimensionless 
parameters are 
$$
\mu = \frac{k_{-1}}{k_1s_0}, \qquad \nu = \frac{k_{-1}+k_2}{k_1 s_0}, 
\qquad \varepsilon = \frac{e_0}{s_0} .
$$ 
From (\ref{inits}), the initial conditions for (\ref{diml}) 
are 
\beq\label{ninits} 
x(0) = 1, \qquad y(0) = 0 = z(0).
\eeq 
Typically, the amount of enzyme  is very small compared with the 
concentrations of the other reactants, so the parameter $\varepsilon$ should be small. According to Murray \cite{murray2}, 
the realistic range is $10^{-7} \leq \varepsilon \leq 10^{-2}$.

\begin{figure}\label{pxy} 
\vspace{0.1cm}
\centerline{
\scalebox{0.5}{\includegraphics{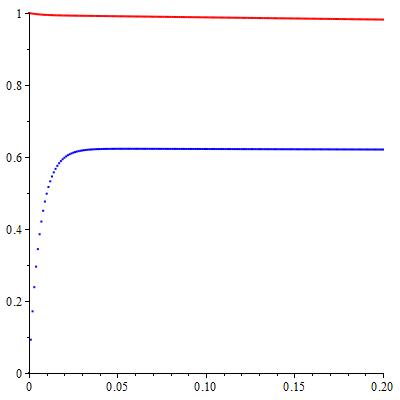}}
}
\caption[ ]{As for Figure 1, with details of $x$ (above) and $y$ (below) against $t$ for $0\leq t\leq 0.2$.}          
\end{figure}

In terms of the dimensionless variables, Kahan's method applied to (\ref{diml}) yields the 
discrete three-dimensional system
\beq \label{ddiml} 
\begin{array}{rcl}
2(\tilde{x}-x)/h & = & -(x+\tilde{x}) +\mu (y +\tilde{y}) + x\tilde{y}+\tilde{x}y, \\ 
2\varepsilon (\tilde{y}-y)/h & = &  x+\tilde{x} -\nu (y +\tilde{y}) - x\tilde{y} - \tilde{x}y, \\ 
2(\tilde{z}-z)/h & = &  (\nu -\mu )(y +\tilde{y}) . 
\end{array}
\eeq 
Using the formula (\ref{phi}) to solve for the shifted variables, the 
resulting map can be written in matrix form as    
$$ 
\left(\begin{array}{c} \tilde{x} \\ \tilde{y} \\ \tilde{z} \end{array} \right) 
= \left(\begin{array}{c} x \\ y \\ z \end{array} \right) + h\,\left(\begin{array}{ccc} 
1+\frac{h}{2}(1-y) & -\frac{h}{2}(\mu + x) & 0 \\ 
-\frac{h}{2\varepsilon}(1-y) & 1+ \frac{h}{2\varepsilon}(\nu +x) & 0 \\ 
0 & -\frac{h}{2}(\nu -\mu ) & 1 
\end{array} \right)^{-1} \left(\begin{array}{c}  
-x + \mu y + xy \\ \varepsilon^{-1} ( x- \nu y - xy ) \\ (\nu - \mu )y \end{array} \right)  
. 
$$ 
In practice, since the equations for $x$ and $y$ decouple from $z$,  the dimensionless concentration 
of the product, it convenient to solve for $x$ and $y$ first, and then the value 
$z_n$ (the value of $z$ after $n$ steps, which approximates $z(nh)$) is given by 
$$ 
z_n = \frac{h}{2}(\nu -\mu ) \sum_{i=0}^{n-1}(y_i + y_{i+1}), 
$$ 
which follows from the third equation in (\ref{ddiml}) and the initial 
conditions (\ref{ninits}).

Some numerical results obtained by applying the discretization (\ref{ddiml}) 
are shown in Figure 1. There are a few things worthy of comment here.  
Note that the value of the time step $h$ has been chosen to be an order of magnitude 
smaller than the small parameter $\varepsilon$. There is also an apparently anomalous feature 
of the left hand plot: the solution curve for $x$ starts with $x=1$ at time zero, as it 
should according to the initial values (\ref{ninits}); but the value of $y$ appears to start at 
$y\approx 0.6$, when it should   start with $y=0$. Subsequently the values of $x$ and $y$ both tend 
asymptotically towards zero; that is to be expected, since the decoupled system for $x$ and $y$, 
given by the first two equations in (\ref{diml}), has a  steady state at $(x,y)=(0,0)$  
which is stable  and  unique, and the results of \cite{roegerstabi}, based on the formula (\ref{mu}), 
imply that the discrete system has the same local stability. From the plot on the right hand side of Figure 1, 
it is seen that the amount of product in the reaction, measured by the variable $z$, tends to an equilibrium value  
(proportional to the area under the graph of $y$).    

The anomalous aspect of Figure 1, namely the initial behaviour of $y$, can be understood by 
looking in detail at the value of $y$ for very small times. This is shown in Figure 2, from which it can be seen 
that in fact the value of $y$ undergoes an early transient phase of very rapid expansion from the 
initial value $y=0$, before reaching a value around 0.6 and then starting to decay.    
The reason for this rapid change is the presence of the parameter $\varepsilon \ll 1$ 
in the equation for $y$. If the left hand side of the second equation in (\ref{diml}) is ignored, then the  
Michaelis-Menten approximation 
$$ 
y\approx \frac{x}{\nu + x} 
$$ 
results, and putting in $\nu = 0.6$ and $x=1$ gives the ``initial'' value $y=0.625$. The initial 
expansion in $y$ takes place over a timescale of the order of $\varepsilon$. A fuller  
understanding of the solution can be achieved by using matched asymptotic expansions, taking 
one type of asymptotic series solution close to time zero (the \textit{inner} solution), and another 
for larger times (the \textit{outer} solution); more details can be found in \cite{murray2}. 
The low resolution of the first part of the plot of $y$ in Figure 2 indicates that a smaller time step $h$ is required 
for the integration at very early times, in order to obtain a more accurate picture of the inner solution.

\section{Birational discretization: Lotka-Volterra systems}

The classic Lotka-Volterra  model for a predator-prey
interaction takes the form 
\beq \label{abcd} 
\begin{array}{rcl} \dot{x} & = & \al x - \be xy , \\ 
\dot{y} & =&  -\ga y +\de xy ,  
\end{array}  
\eeq 
where $x$ and  $y$  denote the sizes of the prey and  predator populations, respectively, and 
the parameters $\al , \be , \ga , \de$ are all positive. With different choices of signs 
for these four parameters, the system can model different types of two-species interaction (e.g. 
competition or mutualism), but once the signs are  fixed then $x,y$ and the time $t$ 
can all be rescaled to obtain a dimensionless model with only a single parameter remaining 
(the other three being set to the value 1).   

In this section we consider discretizations obtained by replacing the $x$, $y$ and $xy$ terms 
appearing in (\ref{abcd}) by expressions of the general  form 
$$ 
x\to a x + \hat{a}\tilde{x}, \qquad y \to    Ay + \hat{A}\tilde{y}, \qquad 
xy \to bxy+c \tilde{x}\tilde{y}+d x\tilde{y} + e \tilde{x}y. 
$$ 
The results in Theorems 1 and 2 below hold true independent of the choice of 
parameters; so henceforth, for convenience, we consider the system with all  parameters set to 1, 
viz  
\beq\label{lv} 
\bear{ccl} 
\dot x & = & x(1-y), \\ 
\dot y & = & y(x-1).  
\eear 
\eeq

The Lotka-Volterra system (\ref{lv}) has a first integral given by 
\beq\label{ham} 
H=\log (xy)-x-y. 
\eeq
This can be viewed as a Hamiltonian function,  
and the system can be given an 
interpretation in terms of a particle moving in one dimension with position $q$ 
and momentum $p$, by 
setting 
$$ 
q= \log x, \qquad p = \log y . 
$$ 
Then expressing the Hamiltonian (\ref{ham}) as a function of $q$ and $p$ gives $H = q+p - e^q - e^p$, 
and the equations (\ref{lv}) can be rewritten in the form 
of a canonical Hamiltonian system: 
$$ 
\dot{q} = \frac{\partial H}{\partial p}, \qquad 
\dot{p} = - \frac{\partial H}{\partial q} . 
$$  
It follows that the flow of (\ref{lv}) is area-preserving in the $(q,p)$ plane, which 
means that (in terms of the original variables $x$ and $y$) the two-form 
\beq\label{om} 
\om  = \frac{1}{xy}\, \rd x \wedge \rd y
\eeq 
is preserved by this flow. 
The trajectories of solutions in the positive quadrant of the $(x,y)$ plane 
are closed curves around the steady state at $(1,1)$, which are level curves 
$H=$constant.

The Lotka-Volterra model was one of the examples originally considered by Kahan 
in his unpublished lectures from 1993. Kahan's method applied to 
(\ref{lv}) yields  
\beq\label{kdlv1} 
\bear{ccl}
(\tilde{x}-x)/{h} & = & \frac{1}{2}\Big( \tilde{x}+x - (\tilde{x}y+x\tilde{y})\Big), \\ 
(\tilde{y}-y)/{h} & = & \frac{1}{2}\Big( \tilde{x}y+x\tilde{y} - (\tilde{y}+y)\Big).  
\eear 
\eeq 
As shown by Sanz-Serna \cite{ss}, the birational map $\varphi : (x,y) \mapsto (\tilde{x},\tilde{y}) $ 
defined by (\ref{kdlv1}) is symplectic, preserving the same two-form (\ref{om}) 
as the original continuous system (\ref{lv}).

Mickens proposed a broad approach to discretization of nonlinear systems, 
with the aim of preserving structural properties of the solutions \cite{mickensbook}.  
In \cite{mickens} he presented a particular discrete predator-prey system,  given by  
\beq\label{mickens} 
\bear{ccl} 
(\tilde{x}-x)/{h} & = & 2x- \tilde{x} - \tilde{x}y, \\ 
(\tilde{y}-y)/{h} & = & - \tilde{y}+2\tilde{x}y-\tilde{x}\tilde{y}.  
\eear 
\eeq 
This gives another explicit birational map of the plane: despite
the fact that the overall  system   is not linear in   $\tilde{x},\tilde{y}$, the first  equation can be solved 
for $\tilde{x}$ and this can be substituted into the 
second equation to obtain  a rational expression for  $\tilde{y}$ in terms of $x$ and $y$. 
(In fact, Mickens uses a general  function $\phi (h)$ in place of $h$ in 
the denominator on the left hand side, but since $\phi (h) = h + o(h)$ this makes 
no difference when $h$ is small.) 
Numerical studies indicate that, for small $h$, 
the discrete system defined by (\ref{mickens}) also has closed orbits around $(1,1)$. 
The map of the $(x,y)$ plane defined by (\ref{mickens})  preserves 
the same symplectic form (\ref{om})  as before, which suggests why its 
stability properties appear to be the same as for the discretization (\ref{kdlv1}). 

For any  system of polynomial differential equations (which need not necessarily be a  quadratic one),
one can try to  implement Mickens' approach 
in the most general way, 
by replacing each monomial $x_ix_j \ldots x_k$ with a linear combination  of all  
possible products of the same variables with/without shifts. 
In \cite{hone} this was done for an example of a cubic vector field, but here 
we illustrate this idea by  applying it to (\ref{lv}), which gives the general discrete system 
\beq\label{genm} 
\bear{ccl} 
(\tilde{x}-x)/{h} & = & ax+ \hat{a} \tilde{x} - (bxy+c \tilde{x}\tilde{y}+d x\tilde{y} + e \tilde{x}y) , \\ 
(\tilde{y}-y)/{h} & = & - Ay -\hat{A}\tilde{y} + (Bxy+C \tilde{x}\tilde{y}+D x\tilde{y} + E \tilde{x}y) .  
\eear 
\eeq 
In order for this to be a first-order method for (\ref{lv}), the coefficients are required to satisfy 
the constraints 
\beq\label{con}  
a+\hat{a} = A + \hat{A} = b+c+d+e = B+C+D+E = 1,  
\eeq 
and it is further assumed that they are all independent of $h$.  
Thus from the start we can say that (in addition to 
the time step $h$) 
the discrete system (\ref{genm}) depends on 8 constant parameters 
$a,b,c,d,A,B,C,D$, with $\hat{a},\hat{A},e,E$ being fixed according to (\ref{con}). 
However, for what follows it will be convenient to specify a system of the form (\ref{genm}) 
by a list of 10 parameters, viz  
$$ 
\{ a,b,c,d,e,A,B,C,D,E\},  
$$
where it is understood that the conditions (\ref{con}) hold, so only 8 of these 10 parameters are independent.  
So for example, Kahan's method produces  the system (\ref{kdlv1}), which is specified by 
$\{ \mez ,  \mez , 0,0, \mez ,  \mez ,  \mez ,0,0,  \mez ,  \mez  \}$, 
while Mickens' system (\ref{mickens}) is specified by $\{ 2,0,0,0,1,0,0,-1,0,2\}$. 

\begin{figure}\label{periodic} 
\begin{center}$ 
\begin{array}{cc} 
{\includegraphics[width=2.2in]{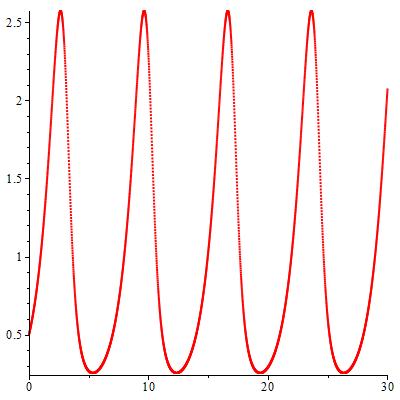}} & 
{\includegraphics[width=2.2in]{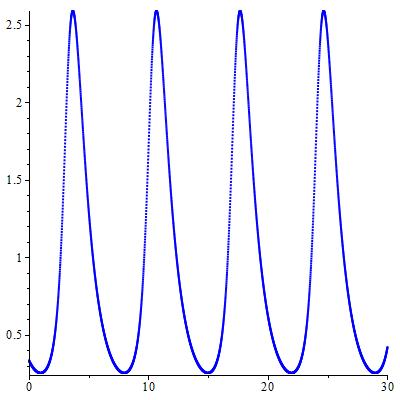}} 
\end{array}$ 
\end{center}
\caption[ ]{Plots of $x$ (left) and $y$ (right)  against $t$  from numerical integration of (\ref{lv}) using the method 
(\ref{genm}) with parameters $\{\mez ,0,\frac{3}{2},-\mez ,0,\mez ,\frac{4}{5},0,\frac{1}{5},0\}$ and $h=0.01$.}  
\end{figure}

The additional requirement that we impose is that the system (\ref{genm}) should 
give a \textit{birational} map, so that $\tilde{x}$ and $\tilde{y}$ can be found explicitly and uniquely 
in terms of $x$ and $y$, and vice versa.   Our main result   can then be stated as follows. 

\begin{theorem}\label{birat}
The system (\ref{genm}) is a birational discrete Lotka-Volterra equation 
if and only if the parameters belong to one of the following cases: 

$$ 
\bear{clcl} 
\mathrm{(i)} & \{ a,0,0,d,e,A,0,0,D,E\}\, & with & d+e=1=D+E; \\ 

\mathrm{(ii)} & \{ a,0,0,1,0,A,B,0,D,E\} &  with  & B+D+E=1; \\  

\mathrm{(iii)}   & \{ a,0,0,0,1,A,0,C,D,E\} & with  & C+D+E=1; \\ 

\mathrm{(iv)} & \{ a,b,0,d,e,A,0,0,0,1\} &  with  & b+d+e=1;  \\  

\mathrm{(v)} & \{ a,0,c,d,e,A,0,0,1,0\} & with  & c+d+e=1; \\  

\mathrm{(vi)} & \{ a,0,c,d,0,A,B,0,D,0\} &  with  & c+d=1=B+D; \\  

\mathrm{(vii)} & \{ a,b,0,0,e,A,0,C,0,E\} &  with & b+e=1=C+E.  
\eear 
$$ 
\end{theorem} 
In order to obtain the above list of parameter choices, in the subsequent subsections we proceed to 
present  two different methods for finding birational discretizations, the first of 
which is a general method  that is applicable to arbitrary vector fields, while the second   
is specific to the quadratic nature of the vector field (\ref{lv}). 

\subsection{The elimination method} 

The first method is to perform successive elimination of variables from the system (\ref{genm}), 
and then impose conditions on the polynomials which result. 
Let us begin by choosing to eliminate $\tilde{x}$. In the example at hand, both equations 
in (\ref{genm}) are linear in  $\tilde{x}$, so we can solve either one explicitly for this variable 
and substitute it into the other equation, to get a polynomial relation between the remaining 
variables $x,y,\tilde{y}$. More generally, if the degree of nonlinearity in $\tilde{x}$ were higher, then 
it would be necessary to eliminate this variable from a pair of equations by taking a resultant, 
and for a system in dimension $N$ one should take $N-1$ resultants  to eliminate one of the variables 
and obtain $N-1$ relations between the remaining variables. In this example, the relation found 
by eliminating $\tilde{x}$  has the form of a quadratic in $\tilde{y}$, that is  
\beq\label{quadr} 
p_2(x,y) \tilde{y}^2 + p_1(x,y) \tilde{y} + p_0(x,y) =0, 
\eeq 
where the coefficients $p_j$ are all polynomials in $x$ and $y$. 

One obvious way to obtain  $\tilde{y}$ as a rational function of 
$x$ and $y$ is to require that $p_2$ in (\ref{quadr}) vanishes, in which 
case (provided $p_1\neq 0$)  $\tilde{y}=-p_0(x,y)/p_1(x,y)$. 
To be precise, up to overall scaling we have 
$$ 
p_2 = (cD-dC)h x + c(Ah-1-h) =0 \Rightarrow cD-dC = 0 = c (Ah-1-h), 
$$ 
since all coefficients of  $p_2$ must vanish. The second condition 
above requires that $c=0$ (since $A$ is assumed independent of $h$), 
and then the first condition implies $dC=0$, and 
so we have 
\beq\label{condns} 
c=0 \qquad \mathrm{and} \, \,\mathrm{either} \quad d=0\quad \mathrm{or} \quad  C=0. 
\eeq 
These conditions are sufficient to ensure that   $\tilde{y}$, and hence also  $\tilde{x}$ (which 
can be found in terms of  $x,y$ and $\tilde{y}$ by solving a linear equation), are rational functions 
of $x$ and $y$. Thus we have a rational map $\varphi : (x,y) \mapsto (\tilde{x},\tilde{y}) $ . 

In order for the inverse map  $\varphi^{-1}$ to be rational, we require sufficient conditions for $x$ and 
$y$ to be given as rational functions of   $\tilde{x}$ and $\tilde{y}$. To do this, we choose to eliminate 
$x$ from the system (\ref{genm}), to obtain a quadratic in $y$, of the form 
\beq\label{nquadr} 
 \tilde{p}_2( \tilde{x}, \tilde{y}) {y}^2 +  \tilde{p}_1( \tilde{x}, \tilde{y}){y} +  \tilde{p}_0( \tilde{x}, \tilde{y}) =0, 
\eeq 
where, 
up to scaling,  
\beq \label{pe2t} 
  \tilde{p}_2 = \Big(B (1-c-d) -b(1-C-D) \Big) h  \tilde{x} +b(Ah-1) , 
\eeq 
and we have chosen to remove $e$ and $E$ from the formulae by using the constraints (\ref{con}). 
Now we can  obtain the rational expression 
$y= - \tilde{p}_0( \tilde{x}, \tilde{y})/ \tilde{p}_1( \tilde{x}, \tilde{y})$ 
by requiring that $ \tilde{p}_2 = 0$; from 
(\ref{pe2t}) this gives $b=0$, because $Ah-1=0$ is not possible, and then we 
have 
\beq\label{ncondns} 
b=0 \qquad \mathrm{and} \, \, \mathrm{either} \quad c+d=1\quad \mathrm{or} \quad  B=0. 
\eeq 
If these conditions hold, then $y$ is a rational function of   $\tilde{x}$ and $\tilde{y}$, 
and hence also $x$ is.

Overall we see that requiring both sets of conditions (\ref{condns}) and (\ref{ncondns}) to hold is sufficient 
for the map $\varphi$ to be birational.  This leads to three possibilities, which are the cases (i),(ii) and 
(iii) in Theorem 1. However, other cases are possible if we choose to eliminate the variables in a 
different order. So for instance, eliminating   $\tilde{y}$  first gives a quadratic in  $\tilde{x}$, 
and then requiring the leading coefficient to vanish implies that  $\tilde{x}$ as a rational 
function of  $x$ and $y$, so $\tilde{y}$ is also; and if $y$ is eliminated next and 
the leading coefficient of the resulting quadratic in $x$ is required to vanish, then a 
different set of sufficient conditions for $\varphi$ to be birational are found. These 
conditions lead to three different possibilities, namely case (i) (again), and cases (iv) and (v).  
Similarly, one can perform the elimination of  $\tilde{x}$ together with 
$y$,  or $\tilde{y}$ together with  $x$; each of these options lead to four possibilities, but overall 
only two new cases arise in this way, namely (vi) and (vii).  

An alternative way to obtain cases (i)-(vii) above is  presented in the next subsection,
but before proceeding with this we make some general 
observations about the result. Each of the discrete Lotka-Volterra equations specified in Theorem 1 
depends on four arbitrary parameters (as well as the time step $h$); the choice of parameters 
$a$ and $A$ is arbitrary in every case. Furthermore, although they are independent 
cases, some of them are related to each other by inversion. To see this, observe that the inverse 
of any discretization method (\ref{genm}) is obtained by switching the roles 
of the dependent variables: $(x,y)\leftrightarrow (\tilde{x}, \tilde{y})$. Performing 
this switch results in another method of the same form but with the parameters 
and time step changed as follows: 
$$ 
h\to -h, \quad a \to 1-a, \quad   A \to 1-A, \quad 
b \leftrightarrow c, \quad  B \leftrightarrow C, \quad d \leftrightarrow e, \quad    
D \leftrightarrow E.
$$ 
From this transformation of the coefficients, it is clear that the inverse of 
a case (i) method is another case (i) method, and similarly the other methods 
are related to one another by such a transformation, so that overall the relationships between the different 
cases under inversion can be summarized by 
$$ 
\mathrm{(i)}   \leftrightarrow     \mathrm{(i)}, 
\qquad 
 \mathrm{(ii)}   \leftrightarrow     \mathrm{(iii)}, 
\qquad 
 \mathrm{(iv)}   \leftrightarrow     \mathrm{(v)}, 
\qquad 
\mathrm{(vi)}   \leftrightarrow     \mathrm{(vii)}.  
\qquad 
$$ 

An example of numerical integration of (\ref{lv}) using one of these birational  methods, 
namely a particular instance of case (vi), is shown in Figure 3. Observe that 
the graphs appear to show $x$ and $y$ varying periodically with time $t$. This is 
consistent with the solutions of (\ref{lv}) in the positive $(x,y)$ plane, which 
lie on closed curves   corresponding to  periodic orbits around the centre at $(1,1)$.

 \subsection{Discriminant conditions and symplectic discretizations} 

The elimination method presented above provided sufficient conditions for the 
map $\varphi$ to be birational, and by trying each possible pair of eliminations 
we obtained all seven cases in Theorem 1. However, there is a second way to 
obtain these conditions, which gives a stronger result. In order to be sure that 
these choices of coefficients are necessary and sufficient, we need to consider the 
elimination process in more detail. 

Observe that, after eliminating $\tilde{x}$ 
to find (\ref{quadr}), it is not strictly necessary that $p_2=0$, but merely that 
the quadratic has only rational roots. For example, case (iv) does not arise 
by setting $p_2=0=\tilde{p}_2$, and yet the quadratic (\ref{quadr}) must have one 
rational root corresponding to the formula for $\tilde{y}$ obtained from 
$\varphi$ in this case. One possibility would be to set $p_0=0$, 
so that the quadratic factorizes as 
$$ 
\tilde{y}\, \Big(p_2(x,y)\tilde{y}+ p_1(x,y)\Big)=0, 
$$ 
with a spurious root $\tilde{y}=0$ that can be neglected. However, we can check 
that setting  $p_0=0$ and then, after eliminating $x$, either 
 $\tilde{p}_2=0$ or  $\tilde{p}_0=0$ does not lead to any consistent 
solutions. Thus we should consider the more general possibility that 
(\ref{quadr}) has two rational roots, one of which is spurious, while the 
other correpsonds to the rational formula for $\tilde{y}$ which is one 
component of the map $\varphi$. This possibility arises if and only if the 
discriminant 
$$ 
\Delta(x,y) = p_1(x,y)^2 - 4p_2(x,y)p_0(x,y)
$$ is a perfect square. Similarly, eliminating $x$ gives  the quadratic 
(\ref{nquadr})  
in $y$, and the roots are rational if and only if the corresponding 
discriminant  $\tilde{\Delta} (\tilde{x},\tilde{y}) $ is a perfect square. 
Imposing this condition on the two discriminants 
$\Delta$ and $\tilde{\Delta}$ yields a set of algebraic conditions on the parameters 
in (\ref{genm}) (omitted here for brevity); then only the cases (i)-(vii) listed in Theorem 1 are possible,   
and this completes the proof of the theorem without needing to consider 
any other eliminations.

\begin{figure}\label{decay} 
\begin{center}$ 
\begin{array}{cc} 
{\includegraphics[width=2.2in]{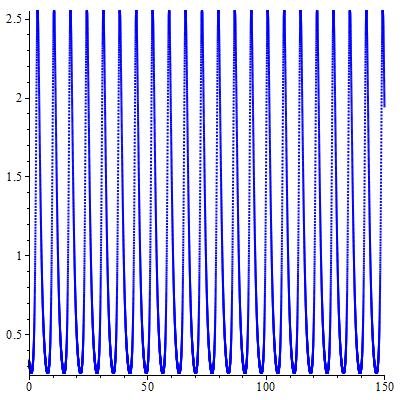}} & 
{\includegraphics[width=2.2in]{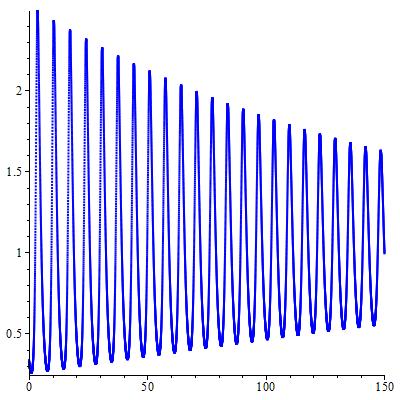}} 
\end{array}$ 
\end{center}
\caption[ ]{Plots of $y$ against $t$  obtained from the method 
(\ref{genm}) with parameters $\{\mez ,\frac{3}{2},0,d ,e,\mez ,0,0,0,1\}$, $e=-d-\mez$ and $h=0.01$, 
for $d=0$ (left) and $d=1$ (right).}  
\end{figure}  

As well as being birational, we should like the map $\varphi$ to preserve the qualitative 
features of the solutions of the continuous system (\ref{lv}). In particular, for the continuous predator-prey 
model, the steady state 
at $ (x,y)=(1,1)$ is a centre in the  phase plane for (\ref{lv}). The form of the 
discretization (\ref{genm}) guarantees that it has the same steady states, but is not enough to 
ensure that the local stability properties are the same. For this model,  Kahan's method sends 
the imaginary eigenvalues $\la = \pm i$ of ${\bf f}'(1,1)$ (the Jacobian of 
(\ref{lv}) at (1,1), that is) to eigenvalues $\mu (h)$ of $\varphi '(1,1)$ which lie on the unit 
circle, as can be seen directly from the formula (\ref{mu}), or less directly by noting that 
(as shown in \cite{ss}) the map 
$\varphi$ is symplectic in this case, so its steady states can only be of centre or saddle type. 
However, a non-symplectic method need not preserve local stability. 
Indeed, Figure  4 
shows a comparison between two different methods of type (iv), where the 
numerical integration is performed over a relatively long time compared with Figure 3: 
in the left hand plot of $y$, the periodic 
oscillations persist, while on the right hand side the oscillations decay towards the 
value $y=1$; the first method is symplectic, while the second is not.

In order for the map $\varphi$ defined by (\ref{genm}) to preserve the symplectic form (\ref{om}), its Jacobian 
$\varphi '$ must satisfy
$$ 
\det \varphi ' = \frac{\tilde{x}\tilde{y}}{xy}.  
$$   
Roeger presented a method to obtain sufficient conditions for this to hold
in \cite{roeger2}. In terms of the parameters in (\ref{genm}), Roeger's method 
leads to the conditions 
\beq\label{sympcon} 
dE-De=cC=dC=cE=bB=bD=eB=0.   
\eeq 
It turns out that all the maps obtained from these conditions are birational. 
Rather than applying the conditions (\ref{sympcon})  
to the general form of (\ref{genm}) and denumerating the possibilities, 
we can instead take the seven cases from Theorem 1 
and calculate the Jacobian in each case, which shows that these are 
actually necessary and sufficient 
conditions for a birational symplectic discretization. 
\begin{theorem}\label{symp}  
The system (\ref{genm}) is a birational symplectic discrete Lotka-Volterra equation, 
preserving the sympletic form (\ref{om}),  
if and only if the parameters belong to one of the following cases: 
$$ 
\bear{clcl} 
\mathrm{(I)} & \{ a,b,0,0,e,A,0,C,0,E\} \, &  with & b+e=1=C+E; \\ 
\mathrm{(II)} & \{ a,0,0,d,e,A,0,0,d,e\}\, & with & d+e=1; \\ 
\mathrm{(III)} & \{ a,0,c,d,0,A,B,0,D,0\} &  with  & c+d=1=B+D. 
\eear 
$$
\end{theorem} 
We have labelled the three symplectic cases (I),(II),(III) in accordance with the result 
stated on p.944 of \cite{roeger2}. To compare with Theorem 1, note that cases (vi) and (vii) 
are symplectic for all choices of  parameter values, and coincide with cases (III) and (I), respectively. 
The method of type (i) is symplectic if and only if $d=D$ (which, due to (\ref{con}), implies $e=E$), 
in which case it reduces to case (II) of Theorem 2. 
Cases (ii)-(v) are not symplectic in general, but if suitable restrictions are made on the 
parameters then they coincide with particular instances of cases (I) or (III). For example, 
both methods used in Figure 4 are of type (iv), but only the one on the left is symplectic, 
corresponding to case (I) with $b=\frac{3}{2}$, $C=1$.    

\section{Conclusions}

The search for discretizations which preserve structural properties of differential 
equations is a fundamental part of numerical analysis \cite{hairer}. 
Non-standard discretizations (including Kahan's method) have been used effectively for 
biological models \cite{elaydiproc,jang,roeger} and in physics (see the review \cite{PPS2}). 
However, the structural features of these methods are still not fully understood.  

In this paper, we have proposed a systematic way to find 
non-standard  discretizations which are \textit{birational}. 
The advantage of birationality is that both the method and its inverse are 
explicit, so the system can be integrated forwards or backwards in time. 
The elimination method presented in subsection 3.1 is applicable 
to arbitary polynomial vector fields; with minor modifications it could 
also be applied to rational vector fields.   
We conjecture that every polynomial (or rational) vector field should admit a birational   
discretization.

In \cite{hone},   
we presented the birational map defined by 
\beq\label{phih}
\varphi_h: \quad \frac{\tilde{x}-x}{h} = a - \frac{1}{2}( \tilde{x}+x)+ x\tilde{x}\tilde{y}, 
\quad \frac{\tilde{y}-y}{h} =b-x^2\tilde{y}, 
\eeq 
together with its inverse. Both $\varphi_h$ and 
$\varphi_{-h}^{-1}$ result from applying the  elimination method to 
Schnakenberg's cubic system,  
\beq\label{trim} 
\dot{x}=a-x+x^2y, \qquad 
\dot{y}=b-x^2y, 
\eeq 
which arises from a trimolecular reaction,  
in contrast to the systems considered with Kahan's method in section 2. 
The system (\ref{trim}) has a Hopf bifurcation, producing a 
limit cycle for suitable values of $a$ and $b$, and numerical 
and analytical results for the map (\ref{phih}) show that these 
features are preserved by the discretization. 
Further details concerning  (\ref{phih}) and its derivation 
will be given elsewhere.

\noindent {\bf Acknowledgments.} 
KT's studentship was funded by the EPSRC and the School 
of Mathematics, Statistics \& Actuarial Science, University of Kent. 


\end{document}